\documentclass[12pt]{article}

\usepackage{amsmath,amssymb,amsbsy,amsfonts,amsthm,latexsym,
                  amsopn,amstext,amsxtra,euscript,amscd}

\begin{document}

\newtheorem{theorem}{Theorem}
\newtheorem{lemma}[theorem]{Lemma}
\newtheorem{claim}[theorem]{Claim}
\newtheorem{cor}[theorem]{Corollary}
\newtheorem{prop}[theorem]{Proposition}
\newtheorem{definition}{Definition}
\newtheorem{question}[theorem]{Open Question}

%%%%%%%%%%%%%%%%%%%%%%%%%
% Alphabet calligraphic %
%%%%%%%%%%%%%%%%%%%%%%%%%
\def\cA{{\mathcal A}}
\def\cB{{\mathcal B}}
\def\cC{{\mathcal C}}
\def\cD{{\mathcal D}}
\def\cE{{\mathcal E}}
\def\cF{{\mathcal F}}
\def\cG{{\mathcal G}}
\def\cH{{\mathcal H}}
\def\cI{{\mathcal I}}
\def\cJ{{\mathcal J}}
\def\cK{{\mathcal K}}
\def\cL{{\mathcal L}}
\def\cM{{\mathcal M}}
\def\cN{{\mathcal N}}
\def\cO{{\mathcal O}}
\def\cP{{\mathcal P}}
\def\cQ{{\mathcal Q}}
\def\cR{{\mathcal R}}
\def\cS{{\mathcal S}}
\def\cT{{\mathcal T}}
\def\cU{{\mathcal U}}
\def\cV{{\mathcal V}}
\def\cW{{\mathcal W}}
\def\cX{{\mathcal X}}
\def\cY{{\mathcal Y}}
\def\cZ{{\mathcal Z}}

%%%%%%%%%%%%%%%%%%%%%%%
% Alphabet blackboard %
%%%%%%%%%%%%%%%%%%%%%%%
\def\A{{\mathbb A}}
\def\B{{\mathbb B}}
\def\C{{\mathbb C}}
\def\D{{\mathbb D}}
\def\E{{\mathbb E}}
\def\F{{\mathbb F}}
\def\G{{\mathbb G}}
\def\H{{\mathbb H}}
\def\I{{\mathbb I}}
\def\J{{\mathbb J}}
\def\K{{\mathbb K}}
\def\L{{\mathbb L}}
\def\M{{\mathbb M}}
\def\N{{\mathbb N}}
\def\O{{\mathbb O}}
\def\P{{\mathbb P}}
\def\Q{{\mathbb Q}}
\def\R{{\mathbb R}}
\def\S{{\mathbb S}}
\def\T{{\mathbb T}}
\def\U{{\mathbb U}}
\def\V{{\mathbb V}}
\def\W{{\mathbb W}}
\def\X{{\mathbb X}}
\def\Y{{\mathbb Y}}
\def\Z{{\mathbb Z}}

\def\ep{{\mathbf{e}}_p}
\def\em{{\mathbf{e}}_m}

\def\scr{\scriptstyle}
\def\\{\cr}
\def\({\left(}
\def\){\right)}
\def\[{\left[}
\def\]{\right]}
\def\<{\langle}
\def\>{\rangle}
\def\fl#1{\left\lfloor#1\right\rfloor}
\def\rf#1{\left\lceil#1\right\rceil}
\def\le{\leqslant}
\def\ge{\geqslant}
\def\eps{\varepsilon}
\def\mand{\qquad\mbox{and}\qquad}

\def\vec#1{\mathbf{#1}}
\def\inv#1{\overline{#1}}
\def\vol#1{\mathrm{vol}\,{#1}}

\newcommand{\comm}[1]{\marginpar{%
\vskip-\baselineskip %raise the marginpar a bit
\raggedright\footnotesize
\itshape\hrule\smallskip#1\par\smallskip\hrule}}

\def\xxx{\vskip5pt\hrule\vskip5pt}

%%%%%%%%%%%%%%%%%%
%% PAPER BEGINS %%
%%%%%%%%%%%%%%%%%%

\title{\bf  Distribution of Modular Inverses and Multiples   of 
Small Integers and the Sato--Tate Conjecture on Average}

\author{ 
{\sc Igor E. Shparlinski} \\
{Department of Computing, Macquarie University} \\
{Sydney, NSW 2109, Australia} \\
{igor@ics.mq.edu.au}}

\date{\today}
\pagenumbering{arabic}

\maketitle

\begin{abstract}
We show that, for sufficiently large integers $m$ and  $X$, 
for almost all $a =1, \ldots, m$ the ratios $a/x$ and the
products $ax$, where $|x|\le X$, 
are very uniformly distributed in the residue ring modulo  $m$. 
This extends some recent results of Garaev and
Karatsuba. We apply this result to show that on average over $r$ and $s$, ranging over 
relatively short intervals,  the
distribution of Kloosterman 
sums
 $$
K_{r,s}(p) = \sum_{n=1}^{p-1} \exp( 2 \pi i (rn + sn^{-1})/p),
$$
for primes $p\le T$ 
is in accordance with the Sato--Tate conjecture. 
\end{abstract}

\section{Introduction}
\label{sec:intro}

\subsection{Motivation} 

A rather old conjecture asserts that if $m = p$ is prime
then for  any fixed $\eps > 0$ and sufficiently large $p$, 
for every integer $a$ there are  integers $x$ and $y$ with $|x|,|y|\le p^{1/2 + \eps}$ 
and such that 
$a \equiv xy \pmod p$,
see~\cite{Gar1,GarKar1,GarKar2,GarKu} and references therein. 
The question has probably been motivated by the following observation. 
Using the Dirichlet
pigeon-hole principle, one can easily show that
for every integer $a$ there are integers $x$ and $y$ with 
$|x|,|y|\le 2p^{1/2}$  with
$a \equiv y/x \pmod p$. 
 
Unfortunately, this  is known only with 
$|x|,|y| \ge C p^{3/4 }$ for some absolute constant $C>0$, which is due to Garaev~\cite{Gar2}.

On the other hand,
it has been shown in the series of works~\cite{Gar1,GarKar1,GarKar2,GarKu}
that the congruence $a \equiv xy \pmod p$ is solvable for all but $o(m)$ values
of  $a = 1, \ldots, m-1$,  with $x$ and $y$ 
significantly smaller than $m^{3/4}$. 
In particular,  it is shown  by Garaev and Karatsuba~\cite{GarKar2} 
for  $x$ and $y$  in the range $1\le x,y\le m^{1/2}(\log m)^{1 + \eps}$. 
Certainly this result is very sharp. 
Indeed, it has been noticed by Garaev~\cite{Gar1} that   well known estimates for
integers with a  divisor   
in a given interval immediately imply that for any $\varepsilon > 0$ almost all residue 
classes modulo $m$ are not of the form $xy \pmod m$ with 
$1 \le x,y \le m^{1/2}(\log m)^{\kappa -\varepsilon}$ where
$$
\kappa = 1 - \frac{1 + \log \log 2}{\log 2} = 0.08607\ldots\,.
$$
One can also derive from~\cite{FKSY} that for any $\varepsilon>0$
the  inequality
$$
\max\{|x|,|y| \ : \ xy \equiv 1 \pmod m\} \ge m^{1/2} (\log m)^{\kappa/2} (\log\log m)^{3/4 - \varepsilon}
$$
holds:
\begin{itemize}
\item for  all positive integers $m \le M$, except for possibly $o(M)$ of
them,

\item for all prime $m = p \le M$ except for possibly  $o(M/\log M)$ of
them. \end{itemize} 

Similar questions about the ratios $x/y$,  
have also been studied,  see~\cite{Gar1,GarKar2,Sem1}. 

\subsection{Our results}

It is clear that these problems are special cases of more 
general questions about   the distribution in small intervals 
of residues modulo $m$ of  ratios $a/x$ and products $ax$, 
where $|x| \le X$. In fact here we consider this more $x$ 
from more general sets $\cX \subseteq [-X, X]$.

Accordingly, 
for integers $a$, $m$,   $Y$ and $Z$ and a set of integers $\cX$,  we  denote  
\begin{eqnarray*}
M_{a,m}(\cX;Y,Z)& = &\# \{ x  \in \cX\ :\ a/x  \equiv y  \pmod m,\\
& & \qquad \qquad \qquad  \qquad
\gcd(x,m)=1, \ y \in  [Z+1,Z+Y]\},\\
N_{a,m}(\cX;Y,Z)& = &\# \{ x  \in \cX \ :\ ax \equiv y \pmod m,\\
& & \qquad \qquad \qquad  \qquad \qquad \qquad
  y \in  [Z+1,Z+Y]\}
\end{eqnarray*} 
where the inversion is always taken  modulo $m$.

We note that although in general the
behaviour  of $N_{a,m}(\cX;Y,Z)$ is similar to the behaviour 
of $M_{a,m}(\cX;Y,Z)$, there are some substantial differences. For example,
if $\cX = \{x \in \Z\ : \ |x|\le X\}$ for some $X \ge 1$, then 
$N_{a,m}(\cX;X,0) = 0$ for all integer $a$ with $m - m/X - 1 < a \le  m-1$,
see the argument in~\cite[Section~4]{Gar1}. 
It is also interesting to remark that  the question of asymptotic behaviour 
of  $N_{a,m}(\cX;Y,Z)$  has some applications to the discrete 
logarithm problem, see~\cite{Sem2}.

Here we extend some of the  results of Garaev and Karatsuba~\cite{GarKar2} 
and show that if  $X,Y \ge m^{1/2 +
\eps}$ and $\cX$ is a sufficiently massive subset of the interval $[-X,X]$,   
then  $M_{a,m}(\cX;Y,Z)$ and  
$N_{a,m}(\cX;Y,Z)$ are  close to their expected average values for all but $o(m)$ values
of  $a  =1, \ldots, m$.

It seems that the method of Garaev and Karatsuba~\cite{GarKar2}  
is not suitable for obtaining results of this kind. So  we use a different approach  
which is somewhat similar to that used in the proof of~\cite[Theorem~1]{BHS}.

Finally we note that one can also obtain analogous results for 
\begin{eqnarray*}
 N_{a,m}^*(\cX;Y,Z)& = &\# \{ x  \in\cX\ :\ ax \equiv y \pmod m,\\
 & & \qquad \qquad \qquad  \qquad
\gcd(x,m)=1, \ y \in  [Z+1,Z+Y]\} 
\end{eqnarray*} 
and several other similar quantities. 

\subsection{Applications}

For  integers $r$ and $s$ and a   prime $p$,  we consider Kloosterman 
sums
$$
K_{r,s}(p) = \sum_{n=1}^{p-1} \ep(rn + sn^{-1})  
$$
where as before the inversion is taken modulo $p$. 
We note that for the complex conjugated sum we have
$$
 \inv{K_{r,s}(p)}  =  K_{-r,-s}(p) =  K_{r,s}(p)
$$
thus $K_{r,s}(p)$ is real. 

Since accordingly to the Weil bound, see~\cite{IwKow,Katz,KatzSar,LN}, we have
$$
\left| K_{r,s}(p)\right| \le 2 \sqrt{p}, \qquad \gcd(r,s,p) =1, 
$$
we can now define the angles $\psi_{r,s}(p)$ 
by the relations
$$
K_{r,s}(p) = 2\sqrt{p} \cos \psi_{r,s}(p)\mand 0 \le \psi_{r,s}(p) \le \pi.
$$

The famous {\it  Sato--Tate\/} conjecture asserts that, 
for any fixed non-zero integers $r$ and $s$, the angles
$\psi_{r,s}(p)$ are distributed accordingly to  the
{\it  Sato--Tate  density\/}
$$
 \mu_{ST}(\alpha,\beta) = \frac{2}{\pi}\int_\alpha^\beta  \sin^2 \gamma\, d \gamma,
$$ 
see~\cite[Section~21.2]{IwKow}.
That is, if $\pi_{r,s}(\alpha, \beta; T)$ denotes the number 
of primes $p\le T$ with $\alpha \le \psi_{r,s}(p) \le \beta$, 
where, as usual $\pi(T)$ denotes the total number of primes $p \le T$, 
the Sato--Tate  conjecture predicts that 
\begin{equation}
\pi_{r,s}(\alpha, \beta; T) \sim  \mu_{ST}(\alpha,\beta)  \pi(T), \qquad T \to \infty,
\end{equation}
for all fixed real
$0 \le \alpha < \beta \le \pi$, see~\cite[Section~21.2]{IwKow}. 
It is also known that if $p$ is sufficiently large and $r$ and $s$ run
independently through $\F_p^*$ then the distribution of $\psi_{r,s}(p)$ is
accordance with the Sato--Tate conjecture, see~\cite[Theorem~21.7]{IwKow}. 
An explicit quantitative bound on the discrepancy between the distribution 
of $\psi_{r,s}(p)$, $r,s \in \F_p^*$ and the Sato--Tate distribution
is given by Niederreiter~\cite{Nied}. 
Various 
modifications and generalisations of this conjecture are given by 
Katz and Sarnak~\cite{Katz,KatzSar}. 
Despite a series of significant efforts towards this 
conjecture, it remains open, for example, 
see~\cite{Adolph,ChaLi,FoMic,FMRS,Katz,KatzSar,Lau,Nied} and references therein.

Here, combining our bounds of $M_{a,m}(\cX;Y,Z)$ with a  
result of Niederreiter~\cite{Nied}, we show that on average over $r$ and $s$, ranging 
over relatively short intervals $|r| \le R$, $|s| \le S$, 
the Sato-Tate conjecture 
holds on average and the sum
$$
\Pi_{\alpha, \beta}(R,S,T) = \frac{1}{4RS}\sum_{0 < |r| \le R} \sum_{0< |s| \le S}
\pi_{r,s}(\alpha, \beta; T)  
$$
satisfies
$$
\Pi_{\alpha, \beta}(R,S,T) \sim  \mu_{ST}(\alpha,\beta) \pi(T). 
$$
Furthermore, over a larger intervals, we also estimate the dispersion
$$
\Delta_{\alpha, \beta}(R,S,T)  = \frac{1}{4RS}\sum_{0 < |r| \le R} \sum_{0< |s| \le S}
\(\pi_{r,s}(\alpha, \beta; T)  - \mu_{ST}(\alpha,\beta) \pi(T)\)^2.
$$

We recall that 
Fouvry and Murty~\cite{FoMu} have the {\it  Lang--Trotter conjecture\/}
on average over $|r|\le R$ and $|s| \le S$  
for the family of elliptic
curves 
$\E_{r,s}$  given  by the {\it  affine Weierstra\ss\
equation\/}:
$$ 
\E_{r,s}~:~U^2 = V^3 + rV + s. 
$$
Several more interesting questions on elliptic curves have 
been studied ``on average'' for similar families of curves 
in~\cite{AkDavJur,Baier,BBIJ,DavPapp1,DavPapp2,Gek,James,JamYu}. 

However, we note that technical details of our 
approach are different
from that of Fouvry and Murty~\cite{FoMu}. For example, 
their result is nontrivial only if  
$$
R S\ge T^{3/2 + \eps} \mand \min\{R, S\} \ge
T^{1/2 + \eps}
$$
for some fixed $\eps> 0$. The technique of~\cite{FoMu} 
can also be applied to getting an asymptotic formula for
 $\Pi_{\alpha, \beta}(R,S,T)$ for the same range of parameters
$R$, $S$ and $T$. Apparently it can also be applied to  
$\Delta_{\alpha, \beta}(R,S,T)$ but certainly in an even narrower range
of parameters. 
On the other hand,  our results  for $\Pi_{\alpha, \beta}(R,S,T)$ and 
$\Delta_{\alpha, \beta}(R,S,T)$ are nontrivial for
\begin{equation}
\label{eq:new threshold 1}
RS \ge T^{1  + \eps} 
\end{equation}   
and 
\begin{equation}
\label{eq:new threshold 2}
 RS \ge T^{2  + \eps}  
\end{equation} 
respectively.  

%%We note that when 
%%$p$ is fixed while $r$  and $s$ run independently over 
%%elements of a finite field of $p$ elements,  the Sato-Tate conjecture
%%for elliptic curves $\# \E_{r,s}(F_p)$  has been established
%%by Birch~\cite{Birch}.  In the sequel,  many generalisations
%%have been obtained by Katz and Sarnak, see~\cite{Katz,KatzSar}.
%% 

 \subsection{Notation}
 
 Throughout the paper, any implied constants in symbols $O$ and $\ll$
 may occasionally depend, where obvious, on the real positive
 parameter  $\eps$ and are absolute otherwise. We recall
 that the notations $U \ll V$ and  $U = O(V)$ are both equivalent to
the statement that $|U| \le c V$ holds with some constant $c> 0$.

We use $p$, with or without a subscript, 
to denote a prime number and use $m$ to denote a positive integer. 

Finally, as usual,  $\varphi(m)$ denotes the Euler function of $m$. 

 \subsection{Acknowledgements}  

The author wishes to thank Moubariz Garaev 
for many  useful discussions. 

This work was supported in part by ARC grant DP0556431.

\section{Congruences}
\subsection{Inverses}

We start with the estimate of the average deviation 
between $M_{a,m}(\cX;Y,Z)$ and its expected value 
taken over $a =1 \ldots, m$.   
If the set $\cX \subseteq [-X,X]$  
is dense enough, for example, 
if $\# \cX \ge X m^{o(1)}$,  this bound is nontrivial for $X,  Y\ge m^{1/2 + \varepsilon}$
for any fixed $\varepsilon > 0$ and sufficiently large $m$.

\begin{theorem}
\label{thm:Invers}  For all positive integers $m $, $X$,  $Y$, an arbitrary integer $Z$
and a set  $\cX \subseteq \{x \in \Z\ : \ |x|\le X\}$, 
$$
\sum_{a=1}^m \left|M_{a,m}(\cX;Y,Z) -   \# \cX_m   \frac{Y}{m}  \right|^2 
\le \# \cX (X+Y)  m^{o(1)}. 
$$
where 
$$
\cX_m = \{ x\in \cX\ : \ \gcd(x,m)=1\}.
$$
\end{theorem}

\begin{proof}  
 We denote 
$$
\em(z)  = \exp(2 \pi i z/m).  
$$
Using the identity
$$
\frac{1}{m}\sum_{-(m-1)/2 \le h \le m/2 } \em(hv)=
\left\{\begin{array}{ll}
1&\quad\text{if $v\equiv 0 \pmod m$,}\\
0&\quad\text{if $v\not\equiv 0 \pmod m$,}
\end{array}
\right.
$$
we write 
\begin{eqnarray*}
\lefteqn{M_{a,m}(\cX;Y,Z) = \sum_{x \in \cX_m} \sum_{y =Z+1}^{Z+Y} 
\frac{1}{m}\sum_{-(m-1)/2 \le h \le m/2 }\em\(h(ax^{-1}-y)\)}\\
& & \qquad =   \frac{1}{m}\sum_{-(m-1)/2 \le h \le m/2 } \sum_{x \in \cX_m} 
\em\(hax^{-1}\)\sum_{y
=Z+1}^{Z+Y} \em(-hy)\\
& & \qquad =  
\frac{1}{m}\sum_{-(m-1)/2 \le h \le m/2 }   \em(-hZ)\sum_{\substack{x=1\\ \gcd(x,m)=1}}^X 
\em\(hax^{-1}\)\sum_{y
=1}^Y \em(-hy).
\end{eqnarray*} 
The term  
corresponding to $h=0$ is 
$$ \frac{1}{m}\sum_{x \in \cX_m} 
  \sum_{y =1}^Y 1 
=  \# \cX_m   \frac{Y}{m}. 
$$
Hence 
$$
M_{a,m}(\cX;Y,Z) -  \# \cX_m   \frac{Y}{m}  \ll \frac{1}{m} E_{a,m}(X,Y)  ,
$$
where 
$$E_{a,m}(X,Y) = \sum_{1 < |h| \le m/2 } \left| \sum_{x \in \cX_m} 
\em\(hax^{-1}\)\right| \left|\sum_{y
=1}^Y \em(-hy)\right|. 
$$
Therefore,  
\begin{equation}
\label{eq:M and E}
\sum_{a =1}^m \left|M_{a,m}(\cX;Y,Z) - \# \cX_m   \frac{Y}{m}  \right|^2 
\le \frac{1}{m^2} \sum_{a =1}^m E_{a,m}(\cX,Y)^2 . 
\end{equation}

We now put $J=\fl{\log (Y/2) }-1$ and define the sets
\begin{eqnarray*}
\cH_0 &=& \left\{h\ | \ 1\le|h|\le \frac{m}{Y} \right\},\\
\cH_j&=&\left\{h\ | \ e^j\frac{m}{Y} <|h|\le e^{j+1}\frac{m}{Y}\right\},
\qquad j =1, \ldots, J,\\
\cH_{J+1}&=&\left\{h\ | \ e^{J+1}\frac{m}{Y} <|h|\le m/2 \right\},
\end{eqnarray*}
(we can certainly assume that $J \ge 1$ since otherwise the 
bound is trivial). 

By the Cauchy inequality we have 
\begin{equation}
\label{eq: E and Ej}
E_{a,m}(\cX,Y)^2  \le (J+2) \sum_{j=0}^{J+1} E_{a,m,j}(\cX,Y)^2, 
\end{equation} 
where
$$
 E_{a,m,j}(\cX,Y) = \sum_{h \in \cH_j} \left| \sum_{x \in \cX_m} 
\em\(hax^{-1}\)\right| \left|\sum_{y
=1}^Y \em(-hy)\right|.
$$
Using the bound
$$
\left|\sum_{y
=1}^Y \em(-hy)\right| = \left|\sum_{y
=1}^Y \em(hy)\right| \ll \min\{Y, m/|h|\}
$$
which holds for any integer $h$ with $0 < |h|\le  m/2$, 
see~\cite[Bound~(8.6)]{IwKow}, we conclude 
that 
$$
\sum_{y =1}^Y \em(-hy)  \ll e^{-j} Y, \qquad j  = 0, \ldots, J+1.
$$
Thus
$$
E_{a,m,j}(\cX,Y) \ll e^{-j} Y \left| \sum_{h \in \cH_j}  \vartheta_h \sum_{x \in \cX_m} 
\em\(hax^{-1}\) \right| ,  \qquad j  = 0, \ldots, J+1, 
$$
for some complex numbers $ \vartheta_h$ with $| \vartheta_h | \le 1$ for $|h| \le m$. 
Therefore, 
\begin{eqnarray*}
\lefteqn{
\sum_{a =1}^m E_{a,m,j}(\cX,Y)^2 \ll  e^{-2j} Y^2 
 \sum_{a=1}^m \left| \sum_{h \in \cH_j}  \vartheta_h \sum_{x \in \cX_m} 
\em\(hax^{-1}\) \right|^2 }\\
&  & \qquad =  e^{-2j} Y^2  \sum_{a=1}^m  \sum_{h_1,h_2 \in \cH_j}  \vartheta_{h_1}\vartheta_{h_2}
\sum_{x_1,x_2 \in \cX_m} 
\em\(a \(h_1x_1^{-1}- h_2x_2^{-1}\)\)\\
&  & \qquad =   e^{-2j} Y^2 \sum_{h_1,h_2 \in \cH_j}  \vartheta_{h_1}\vartheta_{h_2}
\sum_{x_1,x_2 \in \cX_m}  \sum_{a=1}^m  \em\(a \(h_1x_1^{-1}- h_2x_2^{-1}\)\).
\end{eqnarray*} 
Clearly the inner sum vanishes if $h_1x_1^{-1} \not \equiv  h_2x_2^{-1} \pmod m$
and is equal to $m$ otherwise. Therefore
\begin{equation}
\label{eq: sum E_j}
\sum_{a =1}^m E_{a,m,j}(\cX,Y)^2 \ll e^{-2j} Y^2   m T_j ,  
\end{equation}
where $T_j$ is the number of solutions to the congruence
$$
h_1x_2 \equiv  h_2x_1 \pmod m, \qquad 
h_1,h_2 \in \cH_j, \ x_1,x_2 \in \cX_m.
$$
We now see that if $h_1$ and $x_2$ are fixed then $h_2$ and $x_1$ are such 
that their product $s = h_2 x_1 \ll e^j mX/Y$ belongs to a
prescribed residue class modulo $m$. Thus there are at most $O\(e^jX/Y + 1\)$ 
possible values of $s$  and 
for each fixed $s \ll e^j mX/Y$ there are $m^{o(1)}$ values of 
 $h_2$ and $x_1$ with $s = h_2 x_1$,  see~\cite[Section~I.5.2]{Ten}.
Therefore 
$$
T_j \le \# \cX \# \cH_j \(e^jX/Y + 1\) m^{o(1)} = \frac{e^{2j}X  \# \cX m^{1+o(1)}}{Y^2}+ 
 \frac{e^{j} \#\cX m^{1+o(1)}}{Y}
$$
and after substitution into~\eqref{eq: sum E_j}
we get
$$
 \sum_{a =1}^m E_{a,m,j}(\cX,Y)^2 \ll e^{-2j} Y^2   m T_j =  X  \# \cX m^{2+o(1)} + e^{-j}\# \cX Y m^{2+o(1)} .
$$
A combination of this bound with~\eqref{eq: E and Ej} yields the 
inequality
$$
\sum_{a =1}^m E_{a,m}(\cX,Y)^2  \le J^2  X  \# \cX m^{o(1)} +  \# \cX Y m^{2+o(1)} =  \# \cX (X+Y) m^{2+o(1)}.
$$
Finally, recalling~\eqref{eq:M and E}, we conclude the proof.
\end{proof}

\begin{cor}
\label{cor:Inv Interv}  For all positive integers $m$, $X$,  $Y$,  an arbitrary integer $Z$
and the set  $\cX = \{x \in \Z\ : \ |x|\le X\}$ we have 
$$
\sum_{a=1}^m \left|M_{a,m}(\cX;Y,Z) -   2 XY \frac{\varphi(m)}{m^2}  \right|^2 
\le X(X+Y)  m^{o(1)}. 
$$
\end{cor}

\begin{proof}  
 Using the M\"obius inversion formula involving the M\"obius function $\mu(d)$,
see~\cite[Section~1.3]{IwKow} or~\cite[Section~I.2.5]{Ten},  we obtain 
$$ \sum_{\substack{|x| \le X\\ \gcd(x,m)=1}} 1
=\sum_{d | m}\mu (d)\(\frac{2X}{d}+ O(1)\) =  2X \sum_{d | m}\frac{\mu (d)}{d}
+O\( \sum_{d | m}|\mu (d)|\).
$$
Using that 
$$
\sum_{d | m}\frac{\mu (d)}{d} = \frac{\varphi (m)}{m}
$$
see~\cite[Section~I.2.7]{Ten}, and estimating 
$$
\sum_{d | m}|\mu (d)| \le \sum_{d | m} 1 = m^{o(1)}
$$
see~\cite[Section~I.5.2]{Ten}, we derive 
\begin{equation}
\label{eq:Erat} \sum_{\substack{|x| \le X\\ \gcd(x,m)=1}} 1
= 2 X \frac{\varphi (m)}{m}   + O(m^{o(1)}).
\end{equation} 
which  after substitution in Theorem~\ref{thm:Invers} concludes the proof.
\end{proof}

We now immediately derive from Corollary~\ref{cor:Inv Interv}:

\begin{cor}
\label{cor:Almos all M}  For all positive integers $m$, $X$,  $Y$, an arbitrary integer $Z$, 
the set  $\cX = \{x \in \Z\ : \ |x|\le X\}$  and an arbitrary real $\Gamma < 1$,
$$
\left| M_{a,m}(\cX;Y,Z) -  2 XY \frac{\varphi(m)}{m^2}  \right| \ge \Gamma \frac{\varphi(m)}{m^2} XY
$$
for at most $ \Gamma^{-2} Y^{-1}\(X^{-1} + Y^{-1}\) m^{2+ o(1)}$ values of $a =1, \ldots, m$.
\end{cor}

\subsection{Multiples}

We now estimate  the average  deviation 
between $N_{a,m}(\cX;Y,Z)$ and its expected value taken over $a =1, \ldots, m$.  Our arguments
are almost identical to those of Theorem~\ref{thm:Invers}, so we only indicate 
a few places where they differ (mostly only typographically). 
As before, if  $\cX \subseteq [-X,X]$  
is dense enough, for example, 
if $\# \cX \ge X m^{o(1)}$,  this bound is nontrivial for $X,  Y\ge m^{1/2 + \varepsilon}$
for any fixed $\varepsilon > 0$ and sufficiently large $m$.

\begin{theorem}
\label{thm:Mult}For all positive integers $m $, $X$,  $Y$, an arbitrary integer $Z$
and a set  $\cX \subseteq \{x \in \Z\ : \ |x|\le X\}$,  
$$
\sum_{a =1}^m \left|N_{a,m}(\cX;Y,Z) - \# \cX \frac{Y}{m}  \right|^2 
\le \# \cX(X+Y) m^{o(1)}. 
$$
\end{theorem}

\begin{proof}  
As in the proof of Theorem~\ref{thm:Invers}, we write
$$ 
 N_{a,m}(\cX;Y,Z) = \sum_{x  \in \cX} \sum_{y =Z+1}^{Z+Y} 
\frac{1}{m}\sum_{-(m-1)/2 \le h \le m/2 }\em\(h(ax-y)\)
$$
and obtain, instead of~\eqref{eq:M and E},  that 
$$
\sum_{a =1}^m  \left|N_{a,m}(\cX;Y,Z) -   \# \cX \frac{Y}{m}   \right|^2 
\le \frac{1}{m^2} \sum_{a =1}^m  F_{a,m}(\cX,Y)^2 +  Y^2m^{-1 + o(1)}
$$
where
$$F_{a,m}(\cX,Y) = \sum_{1 < |h| \le m/2 } \left|\sum_{x  \in \cX} 
\em\(hax\)\right| \left|\sum_{y
=1}^Y \em(-hy)\right|.
$$
Furthermore, instead of~\eqref{eq: E and Ej}
we obtain 
$$
F_{a,m}(\cX,Y)^2  \le (J+2) \sum_{j=0}^{J+1} F_{a,m,j}(\cX,Y)^2, 
$$ 
where
$$
F_{a,m,j}(\cX,Y) = \sum_{h \in \cH_j} \left|\sum_{x  \in \cX} 
\em\(hax\)\right| \left|\sum_{y
=1}^Y \em(-hy)\right|, 
$$
with the same sets $\cH_j$ as in the proof of Theorem~\ref{thm:Invers}. 
Accordingly, instead of~\eqref{eq: sum E_j} we get 
$$
\sum_{a =1}^m F_{a,m,j}(\cX,Y)^2 \ll e^{-2j}Y^2 m V_j,
$$
where $V_j$ is the number of solutions to the congruence
$$
h_1x_1  \equiv  h_2x_2  \pmod m, \qquad 
h_1,h_2 \in \cH_j, \  x_1 ,  x_2\in   \cX , \ 
\gcd(x_1x_2,m)=1.  
$$
Fixing $h_1$ and $x_1$ and counting the number of possibilities for 
the pair $(h_2,x_2)$, as before,  we obtain
$$
V_j \le  \frac{e^{2j}X \# \cX m^{1+o(1)}}{Y^2} + \frac{e^{j} \# \cX m^{1+o(1)}}{Y},  
$$
which yields the desired result. 
\end{proof}

Using~\eqref{eq:Erat}, we deduce an analogue of Corollary~\ref{cor:Inv Interv}.

\begin{cor}
\label{cor:Mult Interv}  For all positive integers $m$, $X$,  $Y$, an arbitrary integer $Z$
and the set  $\cX = \{x \in \Z\ : \ |x|\le X\}$, 
$$
\sum_{a=1}^m \left|M_{a,m}(\cX;Y,Z) -   2 XY \frac{\varphi(m)}{m^2}  \right|^2 
\le X(X+Y)  m^{o(1)}. 
$$
\end{cor}

We now immediately derive from  Corollary~\ref{cor:Mult Interv}

\begin{cor}
\label{cor:Almos all N}  For all positive integers $m$, $X$,  $Y$, an arbitrary integer $Z$, 
the set  $\cX = \{x \in \Z\ : \ |x|\le X\}$ and an arbitrary real $\Gamma < 1$,
$$
\left| N_{a,m}(\cX;Y,Z) - \frac{2 XY}{m}  \right| \ge \Gamma \frac{XY}{m} ,
$$
for at most $ \Gamma^{-2}  Y^{-1}\(X^{-1} + Y^{-1}\) m^{2+ o(1)}$ values of $a =1, \ldots, m$.
\end{cor}

\section{Distribution of Kloosterman sums}

\subsection{Distribution for a fixed prime}

Let $\cQ_{\alpha,\beta}(R,S,p)$ be the set  
of integers $r$ and $s$  with $|r| \le R$, $|s| \le S$, $\gcd(rs,p)=1$ and 
such that $ \alpha \le \psi_{r,s}(p) \le \beta$.

\begin{theorem}
\label{thm:Kloost}
For all primes $p$ and positive integers $R$ and  $S$, 
$$
\max_{0 \le \alpha < \beta \le \pi} 
\left|\# \cQ_{\alpha,\beta}(R,S,p) - 4\mu_{ST}(\alpha,\beta) RS \right| \ll 
RSp^{-1/4} +  R^{1/2} S^{1/2} p^{1/2+o(1)}.
$$
\end{theorem}

\begin{proof}   Let $\cA_p(\alpha,\beta)$ be the set  
of integers $a$  with $1 \le a \le p-1$  and 
such that $ \alpha \le \psi_{1,a}(p) \le \beta$. 
By the  result of Niederreiter~\cite{Nied}, we have: 
\begin{equation}
\label{eq:Nied Bound}
\max_{0 \le \alpha < \beta < \pi} 
\left|\# \cA_p(\alpha,\beta) - \mu_{ST}(\alpha,\beta) p  \right| \ll p^{3/4}.
\end{equation}
Assume that $R \le S$. Then, using that
$$
K_{r,s}(p) = K_{1,rs}(p), 
$$
and defining the set 
\begin{equation}
\label{eq:Set R}
\cR = \{r \in \Z\ : \ |r| \le R\},
\end{equation}
we write, 
$$
\# \cQ_{\alpha,\beta}(R,S,p) = \sum_{a \in \cA_p(\alpha,\beta)} M_{a,p}(\cR; 2S+1, -S-1) + O\( RS/p \), 
$$
where  the term $ O(RS/p)$ accounts for $r$ and $s$ with $\gcd(rs, p)> 1$. 
Thus the Cauchy inequality and   Theorem~\ref{thm:Invers} yield
\begin{eqnarray*} 
\lefteqn{
\# \cQ_{\alpha,\beta}(R,S,p) - \# \cA_p(\alpha,\beta) \frac{2R(2S+1)}{p}}\\
& & \quad \ll \sum_{a \in \cA_p(\alpha,\beta)}\left| M_{a,p}(\cR; 2S+1, -S-1)- \frac{2R(2S+1)}{p}\right|  +  RS/p \\
& & \quad \ll   \(p\sum_{a=1}^p\left| M_{a,p}(\cR; 2S+1, -S-1)- \frac{2R(2S+1)}{p}\right|^2\)^{1/2}  +  RS/p \\
& & \quad \ll \sqrt{ R(R+S)} p^{1/2+o(1)} +  RS/p . 
\end{eqnarray*}
Using~\eqref{eq:Nied Bound} we see  that for $R \le S$, 
$$
\# \cQ_{\alpha,\beta}(R,S,p) = 4 \mu_{ST}(\alpha,\beta) RS  + O\(RSp^{-1/4} +  R^{1/2} S^{1/2} p^{1/2+o(1)}\)
$$
uniformly over $\alpha$ and $\beta$. 

For  that $R > S$, we write, 
$$
\# \cQ_{\alpha,\beta}(R,S,p) = \sum_{a \in \cA_p(\alpha,\beta)} M_{a^{-1},p}(\cS, 2R+1, -R-1)
$$
where $\cS = \{s \in \Z\ : \ |s| \le S\}$,
and proceed as before.
\end{proof}  

\subsection{Sato--Tate conjecture on average}

We start with an asymptotic formula for $\Pi_{\alpha, \beta}(R,S,T)$

\begin{theorem}
\label{thm: S-T Aver}
For all positive integers $R$,  $S$ and $T$, 
$$
\max_{0 \le \alpha < \beta \le \pi} 
\left|\Pi_{\alpha, \beta}(R,S,T) - \mu_{ST}(\alpha,\beta) \pi(T)\right| \ll   T^{3/4} +  R^{-1/2}
S^{-1/2}T^{3/2+o(1)}$$
\end{theorem}

\begin{proof}   We have
$$
\Pi_{\alpha, \beta}(R,S,T) = \frac{1}{4RS}  \sum_{p \le T}
\# \cQ_{\alpha,\beta}(R,S,p) 
$$
Applying Theorem~\ref{thm:Kloost}, 
after simple calculations we obtain the result.
\end{proof}  

\begin{theorem}
\label{thm: S-T Discr}
For  all positive integers $R$,  $S$ and $T$, 
$$
\max_{0 \le \alpha < \beta \le \pi} 
\Delta_{\alpha, \beta}(R,S,T) \ll  T^{7/4}+ R^{-1/2} S^{-1/2}T^{3+o(1)}$$
\end{theorem}

\begin{proof} For two distinct primes $p_1$ and $p_2$, let $\cA_{p_1p_2}(\alpha,\beta)$ be the set  
of integers $a$  with $1 \le a \le p_1p_2-1$ and such that 
$$ a \equiv a_1 \pmod {p_1} \mand a \equiv a_2 \pmod {p_2}, 
$$
with some  $a_1 \in \cA_{p_1}(\alpha,\beta)$ and $a_2 \in \cA_{p_2}(\alpha,\beta)$. 

Then, with the set $\cR$ given by~\eqref{eq:Set R}, we have 
\begin{eqnarray*}
\lefteqn{\sum_{0 < |r| \le R} \sum_{0< |s| \le S}
\pi_{r,s}(\alpha, \beta; T)^2 }\\ 
& & \quad   =  2 \sum_{p_1< p_2 \le T}  \sum_{a \in \cA_{p_1p_2}(\alpha,\beta)} 
\(M_{a,p_1p_2}(\cR; 2S+1, -S-1) + O\(\frac{RS}{p_1}\)\) \\
& & \qquad\qquad\qquad\qquad\qquad\qquad\qquad\qquad\qquad\qquad\qquad\qquad\qquad + O(RST), 
\end{eqnarray*}
 where  the term $ O(RS/p_1)$ accounts for $r$ and $s$ with $\gcd(rs, p_1p_2)> 1$
and the term $ O(RST)$ accounts for $p_1 = p_2$. 
Therefore, 
\begin{eqnarray*}
\lefteqn{\sum_{0 < |r| \le R} \sum_{0< |s| \le S}
\pi_{r,s}(\alpha, \beta; T)^2 }\\ 
& & \quad   =  2 \sum_{p_1< p_2 \le T}  \sum_{a \in \cA_{p_1p_2}(\alpha,\beta)} 
M_{a,p_1p_2}(\cR; 2S+1, -S-1) + O\(RST^{1 + o(1)}\).  
\end{eqnarray*}

As in the proof  of Theorem~\ref{thm:Kloost}, we derive 
\begin{eqnarray*}
\lefteqn{\sum_{a \in \cA_{p_1p_2}(\alpha,\beta)} 
M_{a,p_1p_2}(\cR; 2S+1, -S-1)}\\
& & \qquad \qquad   = 
 4 \# \cA_{p_1p_2}(\alpha,\beta)  \frac{RS}{p_1p_2} + O\( \sqrt{R S} (p_1p_2)^{1/2+o(1)}\).
\end{eqnarray*}
Thus, using~\eqref{eq:Nied Bound} we obtain 
\begin{eqnarray*}
\lefteqn{\sum_{a \in \cA_{p_1p_2}(\alpha,\beta)} 
M_{a,p_1p_2}(\cR; 2S+1, -S-1)}\\
& & \qquad \qquad   = 4  \mu_{ST}(\alpha,\beta)^2  RS + O\(RS p_1^{-1/4} +
\sqrt{R S}(p_1p_2)^{1/2+o(1)}\).
\end{eqnarray*}
Hence, 
\begin{eqnarray*}
\lefteqn{\sum_{0 < |r| \le R} \sum_{0< |s| \le S}
\pi_{r,s}(\alpha, \beta; T)^2 }\\ 
& & \qquad   =  8   \mu_{ST}(\alpha,\beta)^2  RS \sum_{p_1< p_2 \le T} 1 + O\(RST^{7/4}+ \sqrt{R S}T^{3+o(1)}\)\\
& & \qquad   = 4 \mu_{ST}(\alpha,\beta)^2  RS\pi(T)^2 + O\(RST^{7/4}+ \sqrt{R S}T^{3+o(1)}\).   
\end{eqnarray*}
Combining the above bound with Theorem~\ref{thm: S-T Aver},
we derive the desired result. 
\end{proof} 
 
 Clearly Theorems~\ref{thm: S-T Aver} and~\ref{thm: S-T Discr}
are nontrivial under the conditions~\eqref{eq:new threshold 1}
and~\eqref{eq:new threshold 2}, 
respectively.

%% \section{Remarks}

We also remark that combining~\cite[Lemma~4.4]{FMRS} (taken with $r=1$) together
with the method of~\cite{Nied},  one can prove an asymptotic formula for $\# \cQ_{\alpha,\beta}(1, S, p)$
for $S \ge p^{3/4 + \eps}$ for any fixed $ \eps > 0$. In turn,  this leads to 
an asymptotic formula for $\Pi_{\alpha, \beta}(1,S,T)$ in the same
 range  $S \ge T^{3/4 + \eps}$. However 
it is not clear how to estimate 
$\Delta_{\alpha, \beta}(R,S,T)$ within this approach.

\end{document}